\documentclass[oneside,12pt]{amsart}

\usepackage{amscd,amsmath,latexsym,hyperref}
\usepackage[mathcal]{euscript}

\newtheorem{thm}{Theorem}[section]
\newtheorem{cor}[thm]{Corollary}
\newtheorem{lem}[thm]{Lemma}
\newtheorem{prop}[thm]{Proposition}

\theoremstyle{definition}
\newtheorem{defin}[thm]{Definition}
\theoremstyle{definition}
\newtheorem{exm}[thm]{Example}

\theoremstyle{remark}
\newtheorem*{rem}{Remark}
\newtheorem*{ack}{Acknowledgment}

\newcommand{\A}{{\mathcal A}}
\newcommand{\R}{{\mathbb R}}
\newcommand{\C}{{\mathbb C}}
\newcommand{\Z}{{\mathbb Z}}
\newcommand{\II}{{\mathbb I}}

\newcommand{\T}{{\mathcal T}}

\renewcommand{\ll}{{\ell }}

\DeclareMathOperator{\rep}{rep}
\DeclareMathOperator{\Rep}{Rep}
\DeclareMathOperator{\Hom}{Hom}
\DeclareMathOperator{\Irr}{Irr}

\begin{document}

\title{On representations and $K$-theory of the braid groups}
\author[A.~Adem]{A.~Adem$^{*}$}
\address{Department of Mathematics, University of Wisconsin,
Madison, WI 53706}
\email{\href{mailto:adem@math.wisc.edu}{adem@math.wisc.edu}}
\urladdr{\href{http://www.math.wisc.edu/~adem/}
{www.math.wisc.edu/\~{}adem}}
\author[D.~Cohen]{D.~Cohen$^\dag$}
\address{Department of Mathematics, Louisiana State University,
Baton Rouge, LA 70803}
\email{\href{mailto:cohen@math.lsu.edu}{cohen@math.lsu.edu}}
\urladdr{\href{http://www.math.lsu.edu/~cohen/}
{www.math.lsu.edu/\~{}cohen}}
\author[F.R.~Cohen]{F.R. Cohen$^{*}$}
\address{Department of Mathematics, University of Rochester,
Rochester, New York 14627 }
\email{\href{mailto:cohf@math.rochester.edu}{cohf@math.rochester.edu}}
\urladdr{\href{http://www.math.rochester.edu/~cohf/}
{www.math.rochester.edu/\~{}cohf}}
\thanks{$^{*}$Partially supported by the NSF}
\thanks{$^\dag$Partially supported by grant LEQSF(1999-02)-RD-A-01
from the Louisiana Board of Regents, and by grant MDA904-00-1-0038
from the National Security Agency}



\begin{abstract}
Let $\Gamma$ be the fundamental group of the complement of a
$K(\Gamma,1)$ hyperplane arrangement (such as Artin's pure braid
group) or more generally a \emph{homologically toroidal} group as
defined below.  The subgroup of elements in the complex $K$-theory of
$B\Gamma$ which arises from complex unitary representations of
$\Gamma$ is shown to be trivial.  In the case of real $K$-theory, this
subgroup is an elementary abelian $2$-group, which is characterized
completely in terms of the first two Stiefel-Whitney classes of the
representation.  Furthermore, an orthogonal representation of $\Gamma$
gives rise to a trivial bundle if and only if the representation
factors through the spinor groups.

In addition, quadratic relations in the cohomology algebra of the pure
braid groups which correspond precisely to the Jacobi identity for
certain choices of Poisson algebras are shown to give the existence of
certain homomorphisms from the pure braid group to generalized
Heisenberg groups.  These cohomology relations correspond to
non-trivial Spin representations of the pure braid groups which give
rise to trivial bundles.
\end{abstract}

\maketitle

\section{Introduction}

Given a discrete group $\Gamma$, consider the set of homomorphisms
$\Rep(\Gamma, G)$ where $G$ denotes either the real orthogonal group 
$O(n)$ or the complex unitary group $U(n)$.  There is a natural map 
of sets $\Rep(\Gamma, G) \to [B\Gamma,BG]$, where the target is the 
set of \emph{pointed} homotopy classes of maps from one classifying 
space to the other.

If $G$ is the complex unitary group $U(n)$, there is a stabilization
map $[B\Gamma, BU(n)] \to [B\Gamma, BU]$, which gives rise to a
natural map $\Rep(\Gamma, U(n)) \to [B\Gamma, BU]= KU^0(B\Gamma)$.
Similarly, there is a natural map $\Rep(\Gamma, O(n) ) \to [B\Gamma,
BO] = KO^0(B\Gamma)$.  Note that \emph{reduced} $K$-theory is used
here.

These natural maps motivate the following.

\begin{defin}
The groups $KU^{0}_{\rep}(B\Gamma)$ and $KO^{0}_{\rep}(B\Gamma)$ are
defined to be the subgroups of $KU^{0}(B\Gamma)$ and $KO^{0}(B\Gamma)$
generated by the images, for all $n \ge 1$, of the maps
\[
\Rep(\Gamma,U(n)) \to [B\Gamma,BU]
\qquad \text{and} \qquad
\Rep(\Gamma,O(n)) \to [B\Gamma,BO].
\]
\end{defin}

The purpose of this article is to study these maps, together with
related constructions, for the following class of discrete groups.

\begin{defin}\label{defin: homologically toroidal}
A discrete group $\Gamma$ is said to be \emph{homologically toroidal}
if there is a homomorphism $\mathcal{F}\to \Gamma$ inducing a split
epimorphism in integral homology, where $\mathcal{F}$ is a finite free
product of free abelian groups of finite rank.

Similarly, a topological space $X$ is said to be \emph{homologically
toroidal} if there is a continuous map $\T \to X$ inducing a split
epimorphism in integral homology, where $\T$ is a finite bouquet of
finite dimensional tori.
\end{defin}

Artin's pure braid groups are examples of homologically toroidal
groups.  More generally, this class of groups includes the fundamental
groups of complements of complex hyperplane arrangements which are
aspherical.  Let $\A$ be a hyperplane arrangement, a finite collection
of codimension one affine subspaces in $\C^{\ll}$, with complement
$M(\A)=\C^{\ll}\setminus \bigcup_{H\in \A}H$.  Call $\A$ aspherical,
or a $K(\Gamma,1)$ arrangement, if $M(\A)$ is an Eilenberg-Mac Lane
space of type $K(\Gamma,1)$.  Examples include the orbit configuration
spaces associated to free actions of finite cyclic groups on $\C^*$,
see \cite{DC,Xi}.

Note that if $\Gamma$ is a homologically toroidal group, then the
classifying space $B\Gamma$ is a homologically toroidal space. 
However, the fundamental group of a homologically toroidal space need
not be a homologically toroidal group.  For instance, the complement
of any complex hyperplane arrangement is a homologically toroidal
space, see Proposition \ref{prop:hta}.  But there are arrangements for
which the homology of the fundamental group of the complement, unlike
that of the complement itself, is not finitely generated, see
\cite{Arv,Ra}.

The relationship between representations of homologically toroidal
groups and the real or complex $K$-theory of their classifying spaces
is studied in this article.  In addition to complements of
$K(\Gamma,1)$ arrangements, including the orbit configuration spaces
mentioned previously, some of the results here apply to the
fundamental groups of orbit configuration spaces associated to
elliptic curves, although these groups (if non-trivial) are not
homologically toroidal as their homology groups are not finitely
generated.  These orbit configuration spaces are studied in
\cite{CKX,CX}.

Recall that a representation of $\Gamma$ in $O(n)$ is said to be a
Spin representation if it factors through the natural composite
$Spin(n) \to SO(n) \to O(n)$.  Bundles associated to such
representations are considered in the next two results.

\begin{thm} \label{thm:HTKtheory1}
Let $\Gamma$ be a homologically toroidal group.
\begin{enumerate}
\item If $g: \Gamma \to U(n)$ is a group homomorphism, then the
induced map
\[
B\Gamma \to BU(n) \to BU
\]
is null-homotopic.  Thus $KU^0_{\rep}(B\Gamma)$ is the trivial group,
and the natural complex $n$-plane bundle over $B\Gamma$ obtained from
the representation $g$ is trivial.

\item If $g: \Gamma \to O(n)$ is a group homomorphism, then the
induced map
\[
B\Gamma \to BO(n) \to BO
\]
is essential if and only if at least one of the first two
Stiefel-Whitney classes of~the representation is non-zero.  If the
first two Stiefel-Whitney classes of $g$ vanish, then the natural real
$n$-plane bundle over $B\Gamma$ associated to the representation $g$
is trivial.  Thus this real $n$-plane bundle is trivial if and only if
$g: \Gamma \to O(n)$ lifts to $Spin(n)$.
\end{enumerate}
\end{thm}

\begin{prop} \label{prop:SpinTrivial}
Let $\Gamma=\pi_{1}M(\A)$ be the fundamental group of the complement
of a complex hyperplane arrangement $\A\subset\C^{\ll}$.  If $g:
\Gamma \to O(n)$ is a representation which lifts to $Spin(n)$, then
the vector bundle over $M(\A)$ associated to this representation is
trivial.
\end{prop}

Trivial vector bundles over complements of arrangements arise in a
number of contexts.  For instance, Kohno \cite{K2} develops Vassiliev
invariants of pure braids using a flat connection on such a bundle
over the configuration space $F(\C,n)$ of $n$ ordered points in $\C$,
the complement of the braid arrangement.  These structures also arise
in mathematical physics, for example in work of Drinfel'd \cite{D} and
Kohno \cite{K} in the context of quasi-Hopf algebras and the
Yang-Baxter equations, and in the Drinfel'd-Kohno monodromy theorem
relating the universal $R$-matrix representation of the braid group to
the monodromy of the Knizhnik-Zamolodchikov differential equations,
see \cite{ChP}.  Under certain conditions on the corresponding flat
connection, Schechtman and Varchenko \cite{SV} give solutions of these
equations in terms of generalized hypergeometric integrals defined on
$F(\C,n)$.

The fundamental group of the complement of any $K(\Gamma,1)$
arrangement is homologically toroidal, see Proposition \ref{prop:hta}. 
For such a group $\Gamma$, the group $KO^0_{\rep}(B\Gamma)$ may be
completely computed as follows.

\begin{prop} \label{prop:SWclasses}
Let $\Gamma$ be the fundamental group of the complement of a
$K(\Gamma, 1)$ arrangement and let $\zeta_1$ and $\zeta_2$ be
arbitrary classes in $H^1(\Gamma;\Z/2\Z)$ and $H^2(\Gamma;\Z/2\Z)$. 
Then there is a finite dimensional orthogonal representation of
$\Gamma$ which factors through the abelianization of $\Gamma$ with
first and second Stiefel-Whitney classes given by $\zeta_1$ and
$\zeta_2$ respectively.  Moreover for these groups the Stiefel-Whitney
classes induce an isomorphism
\[
KO^0_{\rep}(B\Gamma)\cong H^1(\Gamma,\Z/2) \oplus H^2(\Gamma,\Z/2).
\]
\end{prop}

This result admits the following specialization to the pure braid
groups.  The abelianization of the pure braid group $P_n$ on $n$
strands is free abelian of rank $m=\binom{n}{2}$.  The mod-$2$
reduction of this abelianization takes values in an elementary abelian
$2$-group which can be regarded as the natural $\Z/2\Z$-torus in
$O(m)$.  In Section \ref{sec:SW&KOrep}, representations of the pure
braid group which factor through a representation of some
$\Z/2\Z$-torus, and which give rise to any values of the first two
Stiefel-Whitney classes are explicitly constructed.  Thus the elements
in the $K$-theory of the classifying space for the pure braid group
which arise from representations are determined completely, and are
very restricted.

The previous proposition admits a natural extension to the level of
classifying spaces which is given in Section
\ref{sec:ClassifyingSpace}.  Namely, the second stage of the Postnikov
tower for $BO$ captures the contribution to $K$-theory arising from
representations of homologically toroidal groups, and there is a
classifying space which captures $KO^{0}_{\rep}(B\Gamma)$ for a
homologically toroidal group $\Gamma$.

On the other hand, the entire $K$-theory of a homologically toroidal
group can be computed easily.  Indeed, if $X$ is a $CW$-complex whose
suspension is homotopy equivalent to a wedge of spheres, then for any
double loop space $Y$, there are group isomorphisms
\[
[X,Y] \to \bigoplus_{q>0}\Hom(H_q(X,\Z),\pi_q Y).
\]
If $\Gamma$ is homologically toroidal and of finite cohomological
dimension, the suspension $\Sigma(B\Gamma)$ has the homotopy type of a
wedge of spheres (see Section 2).  Hence there are isomorphisms
\[
[B\Gamma,BU] \cong \bigoplus_{q>0}\Hom(H_q(\Gamma,\Z),\pi_q BU)
\quad \hbox{and} \quad
[B\Gamma,BO] \cong \bigoplus_{q>0}\Hom(H_q(\Gamma,\Z),\pi_q BO).
\]

Properties of Spin representations of the pure braid group are
summarized in the next theorem.  One noteworthy feature is that this
result supplies Spin representations arising from the structure of the
cohomology algebra for the pure braid group.  The theorem also gives
the (redundant) statement that Spin representations of the pure braid
group give trivial stable vector bundles.

Recall \cite{Arn,C} that the cohomology ring of the configuration
space $F(\mathbb R^k,n)$ of $n$ ordered points in $\mathbb R^k$ is the
quotient of the exterior algebra generated by elements $A_{i,j}$ of
degree $k-1$ for $1 \le j < i \le n$ by the ideal generated by the
relations
\[
A_{i,j}\cdot A_{i,t} - A_{t,j}\cdot A_{i,t} + A_{t,j}\cdot A_{i,j}
\quad \text{for} \quad 1 \le j < t < i \le n.
\]
The relation
$A_{i,j}\cdot A_{i,t} - A_{t,j}\cdot A_{i,t} + A_{t,j}\cdot A_{i,j}$
will be called the ``three term'' relation below, and is the dual of
the Jacobi identity for a certain choice of Poisson algebra \cite{C}.
This relation is used in the next theorem to construct non-trivial
Spin representations of $P_n$.

\begin{thm} \label{thm:Spin}
The three term relation in the cohomology of the pure braid group is
precisely a choice of lifting of the abelianization homomorphism to a
product of generalized Heisenberg groups, and thus corresponds to a
Spin representation of $P_n$.  This representation is non-trivial, but
gives a trivial vector bundle over the configuration space $F(\C,n)$.
\end{thm}

Turning to Artin's full braid group $B_{n}$, the set of irreducible
complex representations, $\Irr(B_n,GL(m,\C))$, was intensively studied
by Formanek and Procesi in \cite{F,FP}.  For $n > 6$, all such
representations in low degrees, $m < n$, are obtained from
specializations of the reduced Burau representation, possibly tensored
with a one-dimensional representation.

The Burau representation is non-trivial, but induces the trivial
element in complex $K$-theory (but not real $K$-theory) when
specialized at $t=1$.  Since $\C^*$ is path-connected, any
specialization of the Burau representation behaves in an analogous
manner.

These two facts may be used to determine the maps in $K$-theory
induced by maps of the stable braid group to $U(n)$.  The next result
is established in Section \ref{sec:FullBraidGroup}.

\begin{prop} \label{prop:FullBraid}
Let $b_n:B_n \to GL(n,\C)$ be given by evaluation of the Burau
representation at a point in $\C^*$.
\begin{enumerate}
\item The induced map $BB_n \to BGL(n,\C)$ is null-homotopic.  Thus
any element in complex $K$-theory induced by the Burau representation
by evaluation at a unit is trivial.  Furthermore, the tensor product
of this specialization of the Burau representation with any other
representation gives a trivial element in complex $K$-theory.

\item  For $m < n$ and $n > 6$, the natural map
$\Irr(B_n,GL(m,\C)) \to KU^0(BB_n)$ is trivial.

\item Evaluation of the Burau representation at $t = 1$ in the real
numbers $\R$ gives $BB_n \to BGL(n,\R)$ which has order $2$, is
injective in mod-2 homology, and has non-vanishing $i$-th
Stiefel-Whitney class for $2i \leq n$.
\end{enumerate}
\end{prop}

One result of Formanek \cite[Lemma 9]{F} states that an irreducible
complex representation of $B_n$ of dimension $n-1$ for $ n > 2$ does
not extend to $B_{n+2}$.  The following is a consequence.

\begin{prop} \label{prop:StableBraid}
Let $\rho:B_{\infty} \to GL(n, \mathbb C)$ be a finite dimensional
complex representation of the stable braid group.  If the restriction
of $\rho$ to some $B_{m}$ is an irreducible representation of
dimension at least $2$, then $\rho$ is the trivial representation.
Thus the induced map in complex $K$-theory is also trivial.
Furthermore, if $\rho$ is any unitary representation, then $\rho$
factors through $\mathbb Z$, the abelianization of $B_{\infty}$, and
the induced element in complex $K$-theory is also trivial.
\end{prop}

Throughout most of this paper, orthogonal or unitary representations
of homologically toroidal groups will be considered.  The final
section contains brief remarks about representations into $GL(n,R)$,
where $R=\mathbb R$ or $\mathbb C$.  Although there are many more
representations, the $K$-theoretic analysis for these more general
representations will not be addressed directly in this paper.

\section{$K$-theory of homologically toroidal groups}
This section addresses the $K$-theory of certain spaces which include
the classifying spaces of homologically toroidal groups.  For brevity,
denote homology with integer coefficients by $H_*X:=H_*(X;\Z)$.

\begin{prop}
Let $X$ be a finite dimensional $CW$-complex such that the suspension
$\Sigma X$ is homotopy equivalent to a bouquet of spheres.  Consider
the group of pointed homotopy classes of maps $[X,\Omega(Y)]$.
\begin{enumerate}
\item The group of pointed homotopy classes of maps $[X,\Omega(Y)]$ is
isomorphic, as a set, to the set
$\bigoplus_{q >0} \Hom(H_qX, \pi_q \Omega(Y))$.

\item The group of pointed homotopy classes of maps $[X, \Omega^2(Y)]$
is isomorphic, as a group, to
$\bigoplus_{q >0} \Hom(H_q X, \pi_q \Omega^2(Y))$.
\end{enumerate}
\end{prop}
\begin{proof}
Notice that the set $[X, \Omega(Y)]$ is isomorphic as a group to
$[\Sigma X, Y]$.  Since $\Sigma(X)$ is homotopy equivalent to a
bouquet of spheres, the underlying set of the group $[\Sigma(X), Y]$
is isomorphic to $\bigoplus_{q >0} \Hom(H_q X, \pi_q \Omega(Y))$.  In
addition, $[X, \Omega^2(Y)]$ and $[\Sigma(X), \Omega(Y)]$ are
isomorphic as groups.
\end{proof}

\begin{cor} \label{cor:KX}
If $X$ is a space which satisfies the above hypotheses, then
\[
KO^0(X) = \bigoplus_{q>0}\Hom(H_q X, \pi_q BO)
\quad \text{and} \quad
KU^0(X)=\bigoplus_{q>0}\Hom(H_q X, \pi_q BU).
\]
\end{cor}
Note that these can be made explicit using Bott periodicity.

\begin{exm}
For the pure braid group $P_{n}$, it is well known that the suspension
of $K(P_n,1)$ is homotopy equivalent to a bouquet of spheres.  Thus
the $K$-theory of this space is obtained from the tensor product of
the integral cohomology for $K(P_n,1)$ with the $K$-theory of a point.
The $K$-theory of $P_{n}$ is then given in terms of integer cohomology
which is torsion free and has Euler-Poincar\'e series
$(1+t)(1+2t)\cdots (1 +[n-1]t)$.  A similar assertion holds for any
complex hyperplane arrangement, since the single suspension of the
complement is homotopy equivalent to a bouquet of spheres, see
\cite{Sc}.
\end{exm}

Now let $\Gamma$ be a homologically toroidal group.  Then, by
definition, there is a homomorphism $w:\mathcal F \to \Gamma$ inducing
a split surjection in integral homology, where $\mathcal F \cong
\coprod_{1\le i\le m}G_i$ is a finite free product of free abelian
groups $G_i$, each of finite rank.  It follows that the homology of
$\Gamma$ with any trivial coefficients is a split summand of that of
$\mathcal F$.  Moreover, for any cohomology theory $E^*$ there is an
induced map $w^*:E^*(B\Gamma) \to E^*(B\mathcal F)$ which is a split
monomorphism.  There is also a natural isomorphism $\bigoplus_{1\leq i
\leq m} \bar E^*(G_i) \to \bar E^*(\mathcal F).$ Note that if
$E^{*}(B\Gamma) \to E^{*}(BG_{i})$ is trivial for every $i$, then
$E^*(B\Gamma)$ must be trivial.

If $\Gamma$ is of finite cohomological dimension, then the classifying
space $B\Gamma$ has the homotopy type of a finite dimensional
$CW$-complex.  Thus, to apply the above corollary to compute the
$K$-theory of a homologically toroidal group, it suffices to establish
the following.

\begin{lem} \label{lem:bouquet}
If $\Gamma$ is a homologically toroidal group, then the suspension of
the classifying space $B\Gamma$ has the homotopy type of a bouquet of
spheres.
\end{lem}
\begin{proof}
Since $\Gamma$ is homologically toroidal, the homology of $\Gamma$ is
a split summand of the homology of
$\mathcal F = \coprod_{1\leq i \leq m} G_i$, where each $G_{i}$ is
free abelian of finite rank.  The suspension $\Sigma(B\Gamma)$ is a
retract of $\Sigma(B\mathcal F)$, which is a finite bouquet of
spheres.  So, in each degree, the homology of $\Sigma(B\Gamma)$ is a
finitely generated free abelian group.

This, together with the fact that every element in the homology of
$\Sigma(B\mathcal F)$ is spherical, implies that there is a bouquet of
spheres $\mathcal S$ which maps to the suspension of $B\Gamma$, giving
a homology isomorphism.

The map $\mathcal S \to \Sigma(B\Gamma)$ is a homotopy equivalence
provided that the suspension of $B\Gamma$ has the homotopy type of a
$CW$-complex.  Since $\Gamma$ is a discrete group, the classifying
space is naturally a $CW$-complex.  Thus so is the suspension.
\end{proof}

The above observations yield

\begin{prop}
If $\Gamma$ is a homologically toroidal group of finite cohomological
dimension, then
\[
KU^0(B\Gamma)\cong \bigoplus_{q>0} H^{2q}(\Gamma,\Z).
\]
\end{prop}

An analogous calculation can be made for $KO^0(B\Gamma)$.

Next, it is shown that the fundamental group of the complement of a
$K(\Gamma,1)$ arrangement is homologically toroidal; this will be our
main source of examples.  Indeed, all examples of homologically
toroidal groups mentioned above may be realized as fundamental groups
of complements of $K(\Gamma,1)$ arrangements.  Let $\A$ be a
hyperplane arrangement in $\C^{\ll}$, with complement $M(\A)$.  Note
that $M(\A)$ has the homotopy type of a finite dimensional
$CW$-complex.

Recall from Definition \ref{defin: homologically toroidal} that a
topological space is homologically toroidal if it admits a map from a
bouquet of tori which induces a split surjection in integral homology.
The above assertion concerning $K(\Gamma,1)$ arrangements follows from
the fact that the complement of any arrangement is homologically
toroidal.  This result may be derived from (co)homological properties
of the complement of an arrangement known from work of Brieskorn
\cite{Br}, Falk \cite{Fa}, and others (see \cite{OT}).  An alternative
argument is given below.

\begin{prop} \label{prop:hta}
For any complex hyperplane arrangement $\A$, the complement $M(\A)$ is
a homologically toroidal space.
\end{prop}
\begin{proof}
Let $\A$ be an arrangement in $\C^{\ll}$ and, without loss of
generality, assume that $\A$ contains $\ll$ linearly independent
hyperplanes.  The proof is by induction on $\ll$.

In the case $\ll=1$, an arrangement $\A\subset\C$ is a finite
collection of points, and the complement is a bouquet of circles.

For general $\ll$, let $\alpha$ be a homology class in
$H_q(M(\A);\Z)$.  It is enough to show that there are maps
$\beta:(S^1)^q \to M(\A)$ such that $\alpha = \sum_{\beta \in J}
\beta_{*}([T])$, where $[T]$ denotes a choice of fundamental class for
the manifold $T=(S^{1})^{q}$ and $J$ is some (finite) indexing set.

If $q < \ll$, let $W$ be a $q$-dimensional subspace of $\C^{\ll}$ that
is transverse to $\A$.  The intersection $W \cap M(\A)$ may itself be
realized as the complement of a hyperplane arrangement $W \cap \A$ in
$W \cong \C^{q}$, so is homologically toroidal by induction.  Since
$W$ is transverse to $\A$ and is $q$-dimensional, the inclusion $i:W
\cap M(\A) \hookrightarrow M(\A)$ induces an isomorphism
$i_{*}:H_{j}(W \cap M(\A);\Z) \to H_{j}(M(\A);\Z)$ in integral
homology for each $j$, $1 \le j \le q-1$, and a surjection
$i_{*}:H_{q}(W \cap M(\A);\Z) \to H_{q}(M(\A);\Z)$ by a Lefschetz-type
theorem (cf.~\cite{GM}).

Now, as is well known, the homology of the complement of an
arrangement $\A$ is torsion free.  Furthermore, the Betti numbers are
determined by the {\em intersection poset} $L(\A)$, the partially
ordered set of multi-intersections of of elements of $\A$, (typically)
ordered by reverse inclusion, with rank function $L(\A) \to \Z$ given
by codimension, see \cite{OT}.  Since $W$ is transverse to $\A$, the
posets $L(\A)$ and $L(W\cap \A$) are identical through rank $q$.
Consequently, the Betti numbers of the complements $M(\A)$ and $W\cap
M(\A)$ are equal in dimensions $0$ through $q$.  Thus the surjection
$i_{*}:H_{q}(W \cap M(\A);\Z) \to H_{q}(M(\A);\Z)$ is, in fact, an
isomorphism.  This yields the result in case $q <\ll$.

It remains to consider the case $q=\ll$.  Assume first that $\A$ is a
central arrangement in $\C^{\ll}$, that is, $\bigcap_{H\in\A}H \neq
\emptyset$.  It is well known that the complement of such an
arrangement is homeomorphic to a product, $M(\A) \cong \C^{*} \times
M(\text{d}\A)$, where $M(\text{d}\A)$ is the complement of a
``decone'' of $\A$, an arrangement in $\C^{\ll-1}$, see \cite{OT}.
Since $M(\text{d}\A)$ is homologically toroidal by induction, it
follows immediately that $M(\A)$ is as well.

If $\A$ is a non-central arrangement in $\C^{\ll}$, let
$\A_{1},\dots,\A_{k}$ be the central subarrangements of $\A$ which
contain $\ll$ linearly independent hyperplanes.  Then for each $i$,
$1\le i \le k$, the intersection $\bigcap_{H\in\A_{i}} H = z_{i}$ is a
point in $\C^{\ll}$.  Let $B_{i}$ be an open ball of radius $\epsilon$
about $z_{i}$ in $\C^{\ll}$.  For $\epsilon$ sufficiently small, the
intersection $B_{i} \cap M(\A)$ is homeomorphic to $M(\A_{i})$, the
complement of the central subarrangement $\A_{i}$, so is homologically
toroidal.

Finally, it is known that the top homology of $M(\A)$ is isomorphic to
the direct sum
\[
H_{\ll}(M(\A);\Z) \cong \bigoplus_{i=1}^{k} H_{\ll}(B_{i}\cap M(\A);\Z)
\cong \bigoplus_{i=1}^{k} H_{\ll}(M(\A_{i});\Z),
\]
see \cite{OT} or \cite{GM}.  Since $M(\A_{i})$ is homologically
toroidal for each $i$, the result follows.
\end{proof}

In particular, the pure braid group $P_{n}$, the fundamental group of
the complement of the braid arrangement
$\A=\{\ker(z_i-z_j),1\le i <j \le n\}$ in $\C^n$, is homologically
toroidal.

\begin{cor}
Let $\A$ be a complex hyperplane arrangement with complement $M(\A)$.
There are isomorphisms
\[
[M(\A), BU] \to \bigoplus_{2j > 0} H^{2j}(M(\A),\Z).
\]
\end{cor}

\begin{rem}
It is well known that the cohomology of the complement of an
arrangement $\A$ is determined by the combinatorial data recorded in
the intersection poset $L(\A)$, see \cite{OT}.  The above result shows
that the complex $K$-theory of the complement of any complex
hyperplane arrangement is combinatorially determined as well.
Similarly, the real $K$-theory of the complement of an arrangement is
combinatorially determined, see Corollary~\ref{cor:KX}.
\end{rem}

Determining which elements in the $K$-theory of a homologically
toroidal group are ``realized'' by representations hinges upon
deciding which elements in the $K$-theory of a finite sum of integers
arise from representations.  In the case of complex $K$-theory, a
homomorphism from a finite sum of integers to $U(n)$ induces a trivial
map in $K$-theory.

\begin{prop} \label{prop:Ctrivial}
Let $A$ be a free abelian group of finite rank.  Then any homomorphism
$A \to U(n)$ induces a trivial map on complex $K$-theory.
\end{prop}
\begin{proof} Assume that $A$ has rank $q$, and write $A \cong \Z^q$.
The image of a representation of $A$ is abelian.  Since an abelian
subgroup of $U(n)$ is conjugate to a subgroup of the group of diagonal
matrices, any homomorphism $A \to U(n)$ factors through a product of
$n$ maps $ \Z \to S^1$.  So consider a homomorphism $A \to S^1$.
Since the target is an abelian group, this map factors as
\[
A \xrightarrow{\sim} \Z^{q} \to (S^1)^{q} \to S^1.
\]

On the other hand, any homomorphism $\Z \to S^1$ is null-homotopic
after passage to classifying spaces.  This suffices.
\end{proof}

\section{ On flat real bundles over a product of circles}
\label{sec:circles}
This section is certainly well known to experts.  The ``bare-hands''
results here are useful in what follows, and are included for
convenience, as well as completeness.

Recall that any real orthogonal representation of a finitely generated
abelian group, not necessarily finite such as $A_{n}\cong\Z^{n}$, is
conjugate to a Whitney sum of one and/or two-dimensional
representations.  Since commuting unitary transformations are
simultaneously diagonalizable and have eigenvalues of length one, any
finite dimensional unitary representation of a finitely generated
abelian group is a sum of one-dimensional unitary representations.
The analogous statement for finite dimensional orthogonal
representations follows by an extension of scalars argument.

Let $\Theta:A_{n} \to O(2)$ be a two-dimensional representation.  The
Stiefel-Whitney classes as well as associated bundles for such a
representation will be addressed below.

\begin{lem}\label{lem:Twice is Trivial}
If there is an element $e$ of $A_{n}$ for which the determinant of
$\Theta(e)$ is $-1$, then for any choice of generators $e_{i}$ of
$A_{n}$, $\Theta(e_{i})$ is given by either
\begin{equation*} \label{eq:matrix}
{\sf M}(\lambda) = \pm \begin{pmatrix}
  -\cos(\lambda) & \sin(\lambda)  \\
   \hfill\sin(\lambda) & \cos(\lambda)
\end{pmatrix}
\text{ for a fixed real number $\lambda$, or }
\pm \II_{2} = \pm \begin{pmatrix}
 1  &  0  \\
 0  &  1
\end{pmatrix}.
\end{equation*}
Thus the representation $\Theta$ is homotopic through representations
to a Whitney sum of one-dimensional representations.  Furthermore,
twice any one-dimensional representation of $A_n$ gives an $SO(2)$
representation, and thus the associated bundle over $BA_n$ is trivial.
\end{lem}
\begin{proof} If there is an element $e$ in $A_{n}$ for which the
determinant of $\Theta(e)$ is not $1$, then there is a basis element,
say $e_1$, such that $\Theta(e_1)={\sf M}(\lambda)$ for some
$\lambda \in \R$.

Note that ${\sf M}(\lambda)$ has order $2$.  Suppose
${\sf N} \in O(2)$ commutes with ${\sf M}(\lambda)$.  If the
determinant of ${\sf N}$ is $-1$, a calculation reveals that
${\sf N}=\pm {\sf M}(\lambda)$.  If the determinant of ${\sf N}$ is
$+1$, a similar calculation shows that ${\sf N}=\pm\II_{2}$.  Thus,
for all $i$, either $\Theta(e_{i})=\pm {\sf M}(\lambda)$ or
$\Theta(e_{i})=\pm\II_{2}$.

Now define $H:[0,1] \times A_{n} \to O(2)$ by the formula
\[
H(t,e_{i})=
\begin{cases}
\pm {\sf M}(t\lambda) &
\text{if $\Theta(e_{i})=\pm {\sf M}(\lambda)$,}\\
\pm\II_{2}  & \text{if $\Theta(e_{i})=\pm\II_{2}$.}
\end{cases}
\]

Notice that the elements appearing as values of $H(t, e_i)$ all
commute for any fixed value of $t$.  Thus there is a unique extension
of $H$ such that $H(t, -)$ is a group homomorphism.  Consequently,
$\Theta$ is homotopic through homomorphisms to a map $\rho:A_{n} \to
\Z/2\Z \oplus \Z/2\Z$, where $\Z/2\Z \oplus \Z/2\Z$ is the subgroup of
$O(2)$ consisting of diagonal matrices with diagonal entries $\pm 1$.
Hence $B\rho$ is a sum of two one-dimensional representations.  The
first part of the lemma follows.

Finally, observe that twice any one-dimensional representation of
$A_n$ is an $SO(2)$ representation.  Thus the associated bundle is
trivial by Proposition \ref{prop:Ctrivial}.
\end{proof}

\begin{rem}
The homotopy in the above lemma is a special property of free abelian
groups.  For example, define $\Z \to O(2)$ by the matrix $-\II_{2}$.
This representation factors through $\Z/2\Z$.  There is a homotopy
obtained by
\[
H(t)=\begin{pmatrix}
 \hfill \cos(t) & \sin(t) \\
       -\sin(t) & \cos(t)
\end{pmatrix}.
\]
Notice that this is a homotopy of the map $\Z \to O(2)$ through
representations, but is NOT a homotopy through representations of
$\Z/2\Z$.  Namely the matrix $H(t)$ does not have order $2$ if
$\cos(t)\sin(t)\neq 0$.  Thus the null-homotopy of representations
above is a property which depends on the source group being free
abelian in order for the homotopy to be a group homomorphism at each
level $t$.
\end{rem}

To continue to investigate the homotopy class of the map $B\Theta$,
notice that by an argument as above, it suffices to consider those
that are all $-1$'s.  Namely, an $SO(2)$ representation of $A_n$
induces a map which is null-homotopic after passage to classifying
spaces, and taking Whitney sums.  Thus it suffices to assume that each
representation has non-trivial first Stiefel-Whitney class.

Thus consider a sum $\Theta=\bigoplus \theta_{i}$ of surjective
representations
\[
\theta_i: A_{n} \to \Z/2\Z, \quad 1 \leq i \leq q.
\]
Fix a basis $e_1,\dots,e_n$ for $H^1(A_{n};\Z/2\Z)$.
Then, for each $i$,
\[
w_1(\theta_i) = \sum_{1 \leq j \leq n} x_{i,j}e_j
\]
with $x_{i,j}$ equal to either $0$ or $1$, and for each $i$ there is a
$j$ with $x_{i,j} \neq 0$.  Since each $\theta_i$ is a line bundle,
the higher Stiefel-Whitney classes vanish, $w_{k}(\theta_{i})=0$ for
$k>1$.

The Stiefel-Whitney classes of $\Theta$ are readily recorded in terms
of these data.

\begin{lem} \label{lem:SWformulas}
Let $\Theta=\bigoplus \theta_{i}: A_{n} \to O(q)$ be a sum of line
bundles as above.  Then
\[
w_{k}(\Theta)= \sum_{1\le i_{1} < \cdots < i_{k} \le q}
w_{1}(\theta_{i_{1}}) \cdots w_{1}(\theta_{i_{k}}).
\]
In particular,
\[
w_{1}(\Theta)= \sum_{1 \le i \le q} w_{1}(\theta_{i}) =
\sum_{1 \le i \le q}\sum_{1 \leq j \leq n} x_{i,j}e_j,
\]
and
\[
w_{2}(\Theta)= \sum_{1 \le i < j \le q}
w_{1}(\theta_{i})\cdot w_{1}(\theta_{j}) =
\sum_{1 \le i < j \le q} \sum_{1 \le s \neq t \le n}
(x_{i,s}x_{j,t}+x_{i,t}x_{j,s})e_{s} \cdot e_{t}.
\]
\end{lem}

The hypotheses of the vanishing of the total first, and second
Stiefel-Whitney classes of these bundles then gives the next lemma.

\begin{lem} \label{lem:vanish}
Let $\Theta: A_n \to O(q)$ be a sum of line bundles
$\theta_i: A_n \to \Z/2\Z$, $1 \leq i \leq q.$
Assume that
\[
w_1(\Theta) = 0 \quad\text{and}\quad w_2(\Theta) = 0.
\]
If $e_1,\dots,e_n$ is a basis for $H^1(A_n;\Z/2\Z)$ for which
$w_1(\theta_i)=\sum_{1 \leq j \leq n} x_{i,j}e_j$ for each $i$,
$1\le i\le q$,
then the following properties hold:
\begin{enumerate}
\item For each fixed $j$ with $1 \leq j \leq n$,
\[
\sum_{1 \le i \le q} x_{i,j} = 0.
\]
\item For each pair $\{s,t\}$ with $1\le s,t \le n$ and $s \neq t$,
\[
\sum_{1 \le i < j \le q}(x_{i,s}x_{j,t}+x_{i,t}x_{j,s}) = 0.
\]
\end{enumerate}
\end{lem}

The next result characterizes homomorphisms which both arise from
finite dimensional orthogonal representations, and are trivial in real
$K$-theory.

\begin{prop}
Let $\Theta: A_n \to O(q)$ be a sum of line bundles
$\theta_i: A_n \to \Z/2\Z$, $1 \leq i \leq q$.  Then $B\Theta$ in
$[BA_n, BO]$ is trivial if and only if both $w_1(B\Theta)$ and
$w_2(B\Theta)$ vanish.  Furthermore, if $B\Theta$ in $[BA_n, BO]$ is
trivial, then the element $B\Theta$ in $[BA_n, BO(q)]$ is also
trivial.  Thus the bundle associated to $B\Theta$ is a trivial
$n$-plane bundle if and only if the representation $\rho$ lifts to
$Spin(q)$.
\end{prop}

\begin{proof}
Clearly, if $B\Theta$ is trivial, then $w_1(B\Theta)=0$ and
$w_2(B\Theta)=0$.

For the other implication, the map $\Theta: A_n \to O(q)$ factors
through the inclusion of the maximal $\Z/2\Z$-torus given by
$(\Z/2\Z)^q$ by hypothesis.  Thus consider the map $\bar\Theta: A_n
\to (\Z/2\Z )^q$ with image given by an elementary abelian $2$-group
of rank $r$ for some $r \leq q$.  The map $\bar\Theta$ is a sum of
line bundles, $\bar\Theta=\bar\Theta_1 \oplus \bar\Theta_2 \oplus
\cdots \oplus \bar\Theta_r$.

If $r=1$ and $w_1(B\bar\Theta)=0$, then $B\bar\Theta$ is
null-homotopic.  If $r>1$, let $e_1,\dots,e_n$ be a basis for
$H^1(A_n;\Z/2\Z)$, and choose a basis for $(\Z/2\Z)^r$.  Write
\[
w_1(\bar\Theta_i) = \sum_{1\le j \le n} x_{i,j} e_j,\quad 1 \le i \le r,
\]
as before, and let $X=\left(x_{i,j}\right)$ be the corresponding
matrix of coefficients.

By Lemma \ref{lem:vanish}, the vanishing of $w_1(B\bar\Theta)$ implies
that each column of $X$ has an even number of non-zero entries.  After
an appropriate change of basis in $(\Z/2\Z)^r$ if necessary, the
matrix of coefficients is of the form
\[
X= \begin{pmatrix}
0      & \cdots & 0      & x_{1,s} & x_{1,s+1} & \cdots & x_{1,n}\\
0      & \cdots & 0      & x_{2,s} & x_{2,s+1} & \cdots & x_{2,n}\\
0      & \cdots & 0      & 0       & x_{3,s+1} & \cdots & x_{3,n}\\
\vdots &        & \vdots & \vdots  & \vdots    &        & \vdots\\
0      & \cdots & 0      & 0       & x_{r,s+1} & \cdots & x_{r,n}
\end{pmatrix},
\]
with $x_{1,s}=x_{2,s}$ non-zero for some $s$, $1\le s \le n$.  Note
that $w_1(B\bar\Theta)=0$ implies that
$x_{1,t}+x_{2,t}+x_{3,t}+\cdots +x_{r,t}=0$ for each $t$,
$s < t \le n$.

By the above observation and Lemma \ref{lem:SWformulas}, for
$s < t \le n$, the coefficient of $e_s\cdot e_t$ in $w_2(B\bar\Theta)$
is given by
\begin{align*}
\sum_{1 \le i < j \le r}(x_{i,s}x_{j,t}+x_{i,t}x_{j,s})&=
x_{1,s}x_{2,t}+x_{1,t}x_{2,s}+(x_{1,s}+x_{2,s})(x_{3,t}+\cdots+x_{r,t})\\
&=x_{1,s}x_{2,t}+x_{1,t}x_{2,s} + (x_{1,s}+x_{2,s})(x_{1,t}+x_{2,t})\\
&=x_{1,s}x_{1,t}+x_{2,s}x_{2,t}.
\end{align*}

Since $x_{1,s}=x_{2,s} \neq 0$, the vanishing of $w_2(B\bar\Theta)$
implies that $x_{1,t}=x_{2,t}$ for each $t$, $s < t \le n$.  Thus
$\bar\Theta_1=\bar\Theta_2$.  It follows that
$B(\bar\Theta_1 \oplus \bar\Theta_2)$ is null-homotopic, which
establishes the result in the case $r=2$.

For $r>2$, inductively assume that the result holds for $r-1$ with
$r-1 \ge 2$.  Observe that $B\bar\Theta$ is homotopic to
$B(\oplus _{i}\bar\Theta_i)=B(\bar\Theta_1\oplus\bar\Theta_2\oplus\Psi)$,
where $\Psi = \bar\Theta_3 \oplus \cdots \oplus \bar\Theta_r$.
Since $B(\bar\Theta_1 \oplus \bar\Theta_2)$ is null-homotopic,
$B\bar\Theta$ is homotopic to $B\Psi$, and the result follows by
induction.

Notice that the above null-homotopy is supported directly on the
Whitney sum of the bundles above, and does not require stabilization.
Thus, if the first two Stiefel-Whitney classes vanish, the bundle is
trivial.
\end{proof}

\section{Stiefel-Whitney classes and $KO^0_{\rep}(B\Gamma)$}
\label{sec:SW&KOrep}
The purpose of this section is to prove the following theorem, and a
number of related results, including Theorem \ref{thm:HTKtheory1} and
Proposition \ref{prop:SpinTrivial}.

\begin{thm} \label{thm:SW&KOrep}
Let $\Gamma$ be a homologically toroidal group, and let
$\rho: \Gamma \to O(n)$ be a representation with vanishing first two
Stiefel-Whitney classes.  Then $B\rho: B\Gamma \to BO(n)$ and $B\rho:
B\Gamma \to BO$ are null-homotopic.
\end{thm}
\begin{proof}
Let $\rho: \Gamma \to O(n)$ be a representation.  If $\Gamma$ is
homologically toroidal, there is a group
$\mathcal F= \coprod_{1\le i\le m}G_{i}$ with each group $G_{i}$ free
abelian group of finite rank, and a map $w:\mathcal F \to \Gamma$
which induces a split epimorphism on homology (with any trivial
coefficients).  Consequently, the composite
\[
G_{i} \to \Gamma \to O(n) \to O
\]
has vanishing first two Stiefel-Whitney classes for every $i$.  Hence
the induced map $\alpha:G_{i} \to O$ gives the trivial element in
$[BG_{i},BO]$ for every $i$.  It follows that the element $B\rho$ is
trivial in $[B\Gamma,BO]$.

The proof that the map $B\rho:B\Gamma \to BO(n)$ is null-homotopic,
and that the associated bundle is trivial, is given in the next lemma
below.
\end{proof}

\begin{lem} \label{lem:Stable implies Non-stable}
Let $\Gamma$ be a homologically toroidal group, and let $\T$ be a
bouquet of tori, together with a map $g:\T \to K(\Gamma,1)$ which
induces a split surjection in integral homology.  Let $K$ denote the
mapping cone of $g:\T \to K(\Gamma,1)$.  Then, the following hold.
\begin{enumerate}
\item There is a cofibre sequence
\[
\T \to K(\Gamma,1) \to K \to \Sigma(\T) \to \Sigma(K(\Gamma,1))
\to \cdots.
\]

\item The mapping cone $K$ of $g: \T \to K(\Gamma,1)$ is a retract of
$\Sigma(\T)$, and is a co-H-space.

\item The natural map $K(\Gamma,1) \to K$ is null-homotopic,

\item There is a short exact sequence of sets
\[
\to [K,X] \to [K(\Gamma,1),X] \to [\T,X] \to .
\]
Thus the inverse image of the class of the constant map in the set
$[\T,X]$ is the class of the constant map in the set
$[K(\Gamma,1),X]$.  \item If $X = BO(n)$ and the element $f$ in
$[K(\Gamma,1),BO(n)]$ restricts to the trivial element in the set
$[\T,BO(n)]$, then $f$ is trivial.
\end{enumerate}
\end{lem}
\begin{proof}
The mapping cone $K$ of $g$, together with the Barratt-Puppe sequence,
gives a cofibre sequence
\[
\T \to K(\Gamma,1) \to K \to \Sigma(\T) \to \Sigma(K(\Gamma,1))
\to \cdots.
\]
It is shown below that there is a ``section'' for
$\Sigma(\T) \to \Sigma(K(\Gamma,1))$.  That is, $\Sigma(K(\Gamma,1))$
is a retract of $\Sigma(\T)$ via the retraction $\Sigma(\T) \to
\Sigma(K(\Gamma,1))$.

Given any cofibre sequence of path-connected $CW$-complexes
$A \to B \to C$ for which the map $B \to C$ is a retraction, there is
an induced map $A \bigvee C \to B$ (the composite $A \bigvee C \to B
\bigvee B \to B$).  This map induces a homology isomorphism.  Thus
there is a map
\[
K \bigvee \Sigma(K(\Gamma,1)) \to \Sigma(\T)
\]
which is a homology isomorphism.

Notice that the map $\T \to K(\Gamma,1)$ is a surjection on
fundamental groups, and so $\pi_1(K)$ is trivial by the Seifert-Van
Kampen theorem.  Thus the map $K \bigvee \Sigma(K(\Gamma,1)) \to
\Sigma(\T)$ is a homotopy equivalence, and so the ``boundary map'' $K
\to \Sigma(\T)$ is null-homotopic.

For a path-connected $CW$-complex $X$, the Barratt-Puppe sequence
induces an exact sequence of sets (but not necessarily groups)
\[
\cdots \to [K,X] \to [K(\Gamma,1),X] \to [\T,X]\to\{*\}.
\]
Consequently, the inverse image of the class of the constant map
$\T \to X$ is in the image of the map $[K,X] \to [K(\Gamma,1),X]$.
This last map is constant as the map $K(\Gamma,1) \to K$ is
null-homotopic by the above remarks.

If $X = BO(n)$ and $f \in [K(\Gamma,1),BO(n)]$ restricts to a trivial
element in $[\T,BO(n)]$, then $f$ ``comes from'' an element in
$[K,BO(n)]$.  But the map $[K,X] \to [K(\Gamma,1),X]$ is trivial by
the above.  Thus the element $f$ in $[K(\Gamma,1),BO(n)]$ was trivial
to start.

Thus bundles (not just stable bundles) which restrict trivially to
$\T$ are trivial bundles.  The lemma follows.
\end{proof}

\begin{prop} \label{prop:CtrivialPROP}
Let $\Gamma$ be a homologically toroidal group.
\begin{enumerate}
\item Any group homomorphism $\Gamma\to U(n)$ induces a null-homotopic
map $B\Gamma\to BU$, and hence a trivial map in complex $K$-theory.

\item None of the non-trivial elements in $[B\Gamma, BU]$ are induced
by representations of $\Gamma$.

\item Every element in $[B\Gamma, BO]$ which arises as a
representation has order $2$.
\end{enumerate}
\end{prop}
\begin{proof}
If $\mathcal{F}$ is a finite free product of free abelian groups of
finite rank as in the definition of a homologically toroidal group,
then by Proposition~\ref{prop:Ctrivial}, any homomorphism
$\mathcal{F}\to U(n)$ induces a trivial map in $K$-theory.  Given any
epimorphism $\mathcal{F}\to\Gamma$ which induces a split surjection in
homology, together with any homomorphism $\Gamma\to U(n)$, the induced
map $B \Gamma \to BU$ is null-homotopic as observed earlier since
$B\Gamma$ is a split summand in $B\mathcal{F}$ after suspension.

By Lemma~\ref{lem:Twice is Trivial}, twice any orthogonal
representation of $A_n$ gives a null-homotopic map after passage to
classifying spaces as these are sums of $SO(2)$ representations.  Thus
twice any orthogonal representation of $\Gamma$ gives a null-homotopic
map after passage to classifying spaces.
\end{proof}

Observe that Theorem \ref{thm:HTKtheory1} follows from
Theorem \ref{thm:SW&KOrep} and Proposition \ref{prop:CtrivialPROP}.

Now recall that Proposition \ref{prop:SpinTrivial} asserts that any
bundle over the complement $M(\A)$ of a complex hyperplane arrangement
$\A$ which arises from a Spin representation of the fundamental group
$\Gamma=\pi_1 M(\A)$ is a trivial bundle (even though the complement
might not be a $K(\Gamma,1)$ space).
\begin{proof}[Proof of Proposition \ref{prop:SpinTrivial}]
By Proposition \ref{prop:hta}, the complement $M(\A)$ is a
homologically toroidal space.  So there is a bouquet of tori $\T$, and
a map $\T \to M(\A)$ which induces a split epimorphism in integral
homology.  Let $A_n$ denote the fundamental group of one of these
tori.  The associated bundle over any one of these tori is induced by
a representation as the natural map from a torus to $K(\Gamma,1)$ is
homotopic to $B\rho$ for some choice of homomorphism $A_n \to \Gamma$.
Thus if the original bundle is induced by a representation of the
fundamental group, then the associated bundle obtained over the torus
arises as a Spin representation of the fundamental group of the torus.

By Proposition~\ref{prop:CtrivialPROP}, a Spin representation of the
fundamental group of a wedge of tori induces a trivial stable bundle
over the wedge of tori.  By Lemma~\ref{lem:Stable implies Non-stable}
the bundle over $M(\A)$ is trivial.  The proposition follows.
\end{proof}

Next, recall that Proposition \ref{prop:SWclasses} asserts that any
two classes in the first and second cohomology of the fundamental
group of the complement of any $K(\Gamma, 1)$ complex hyperplane
arrangement may be realized as the first and second Stiefel-Whitney
classes of a representation.  This is a consequence of the following.

\begin{prop} \label{prop:SW12}
Let $\A$ be a complex hyperplane arrangement with complement $M(\A)$.
Given cohomology classes $\zeta_i \in H^i(M(\A);\Z/2\Z)$, $i=1,2$,
there is a vector bundle $\xi$ over $M(\A)$ for which the first and
second Stiefel-Whitney classes are $w_1(\xi)=\zeta_1$ and
$w_2(\xi)=\zeta_2$.
\end{prop}
\begin{proof}
Assume that $\A$ is an arrangement of $n$ hyperplanes in $\C^\ll$.
The cohomology of $M(\A)$ is isomorphic to the Orlik-Solomon algebra
of $\A$, which is a quotient of an exterior algebra on $n$ generators,
see \cite{OT}.  This may be realized topologically as follows: The
complement $M(\A)$ is homeomorphic to a slice $W \cap (\C^*)^n$ of the
complex $n$-torus, where $W$ is an appropriate $\ll$-dimensional
affine subspace of $\C^n$.  It is readily checked that the inclusion
$i:W \cap (\C^*)^n \hookrightarrow (\C^*)^n$ induces a surjection in
cohomology.

Now if $\xi$ is a vector bundle over $(\C^*)^n$, then the
Stiefel-Whitney classes of the induced bundle $i^*(\xi)$ over
$M(\A)=W \cap (\C^*)^n$ satisfy $w_j(i^*(\xi)) = i^*(w_j(\xi))$.
Thus, since $i^*$ is surjective, it suffices to prove the proposition
in the case where $\A$ is the arrangement of coordinate hyperplanes in
$\C^n$ and $M(\A)=(\C^*)^n$.  This is a straightforward exercise using
the calculations of Section \ref{sec:circles}.
\end{proof}

Observe that the above vector bundles over $(\C^*)^n$ arise from
representations of the fundamental group $\pi_1((\C^*)^n) = A_n$.  For
any arrangement $\A$ of $n$ hyperplanes, the abelianization of the
fundamental group, $\Gamma=\pi_1M(\A)$, of the complement is
isomorphic to $A_n$.  So if $\A$ is a $K(\Gamma,1)$ arrangement, the
above result shows that any two classes in the first and second
cohomology of $\Gamma$ may be realized as the first and second
Stiefel-Whitney classes of an orthogonal representation which factors
through the abelianization of $\Gamma$.

\begin{exm}
Taking $\A$ to be the braid arrangement, the above results show that
the only possible non-trivial elements in the real $K$-theory of the
pure braid group that arise from representations are torsion elements.
The relevant representations, and associated bundles, will be
constructed explicitly below.

Recall that the abelianization of the pure braid group $P_n$ is free
abelian of rank $m=\binom{n}{2}$.  The abelianization map
$P_n \to \Z^{m}$, followed by mod-$2$ reduction, yields a homomorphism
$P_n \to (\Z/2\Z)^{m}$.  Projection to products of the form
$(\Z/2\Z)^q$ followed by inclusions as $2$-tori in orthogonal groups
or special orthogonal groups give families of useful representations.
Stiefel-Whitney classes of some of these will now be discussed.

For $1\le j < i \le n$, let $A_{i,j} \in H^{1}(P_{n};\Z/2\Z)$ be the
cohomology class dual to the standard generator $\gamma_{i,j}$ which
links strands $i$ and $j$.  Recall that a basis for $H^t(P_n)$ is
given by products
$A_{I,J}=A_{i_1,j_1}A_{i_2,j_2}\cdots A_{i_t,j_t}$ where
$1 \le i_1 < i_2 < \cdots <i_t \le n$ and $j_{m}<i_{m}$ for each~$m$.

Certain choices of maps are listed next.  These maps may be used to
give an explicit proof for the pure braid groups of the existence of
representations which have any fixed choice of first, and second
Stiefel-Whitney class in Proposition \ref{prop:SWclasses}.

Associated to $A_{I,J}\in H^{t}(P_{n};\Z/2\Z)$, define representations
\begin{enumerate}
\item $ \alpha(A_{I,J}):P_n \to O(t)$ to be the composite of the
projection specified by $A_{I,J}: P_n \to (\Z)^t \to (\Z/2\Z)^t$
followed by the inclusion of $(\Z/2\Z)^t$ in $O(t)$; and

\item $\beta(A_{I,J}):P_{n} \to SO(t+1)$ to be the composite of the
projection specified by $A_{I,J}: P_n \to (\Z)^t \to (\Z/2\Z)^t$
followed by the inclusion of $(\Z/2\Z)^t$ in $SO(t+1)$.
\end{enumerate}

For $A_{I,J} = A_{i_1,j_1}A_{i_2,j_2}\cdots A_{i_t,j_t} \in
H^{t}(P_{n};\Z/2\Z)$, the (first two) Stiefel-Whitney classes of the
representations $\alpha(A_{I,J})$ and $\beta(A_{I,J})$ are given by
\[
\begin{aligned}
\alpha(A_{I,J})^*(w_1)&=\sum_{m=1}^{t} A_{i_m,j_m},\qquad
&\alpha(A_{I,J})^{*}(w_{2})
&=\sum_{1 \le m < n \le t} A_{i_m,j_m}A_{i_n,j_n},\\
\beta(A_{I,J})^*(w_1)&=0,
&\beta(A_{I,J})^*(w_2)
&=\sum_{1 \le m < n \le t} A_{i_m,j_m}A_{i_n,j_n}.
\end{aligned}
\]

By choosing appropriate basis elements $A_{I,J}$, and forming Whitney
sums, one can use the representations $\alpha(A_{I,J})$ and
$\beta(A_{I,J})$ to produce a representation for which the first two
Stiefel-Whitney classes are any two classes in $H^{1}(P_{n};\Z/2\Z)$
and $H^{2}(P_{n};\Z/2\Z)$.
\end{exm}

A similar explicit construction can be carried out for an arbitrary
$K(\Gamma,1)$ arrangement.

\section{Pure braid group representations in generalized Heisenberg
groups}
The point of this section is to show how the relations in the
cohomology algebra for the pure braid groups are equivalent to the
existence of certain choices of homomorphisms to groups which
themselves are generalizations of the classical Heisenberg group, and
which arise from Riemann surfaces.  The main point is that the
standard quadratic relations in the cohomology of the pure braid
groups correspond to liftings.  These points are described below, and
are variations of the representations used in the previous section.

First define the generalized Heisenberg group $\pi_g$ as the central
extension
\[
1 \to \Z \to \pi_g \to \Z^{2g} \to 1
\]
with characteristic class $\chi_g$ given by the cup product form for
the cohomology ring of a closed orientable surface of genus $g$.  This
characteristic class is given explicitly by
$\chi_g = \sum_{1\leq i\leq g} x_{i} y_{i}$, where
$x_i,y_i$, $1 \le i \le g$, generate the integral cohomology ring of
$\Z^{2g}$, which is an exterior algebra on these generators.  The
classifying space $B\pi_g$ is a $2$-stage Moore-Postnikov tower.  Thus
$\pi_g$ is a torsion free nilpotent group of nilpotence class $2$.

Notice that there is a ``mod-$n$" reduction of $\pi_g$, say $\bar
\pi_g$, obtained by both replacing $\mathbb Z$ by $\mathbb Z/n\mathbb
Z$, and the characteristic class $\chi_g$ by the naturally associated
mod-$n$ reduction.  Note that there is a morphism of group extensions:

\[
\begin{CD}
\Z @>{}>>   \pi_g @>{}>>    \Z^{2g} \\
  @VVV            @VVV                    @VVV \\
\Z/ n\Z @>{}>>   \bar \pi_g @>{}>>  (\Z/n\Z)^{2g}
\end{CD}
\] \vskip .2in

The group $ \bar \pi_g $ is an example of the tensor product $G
\otimes_{\mathbb Z} R$ for a discrete group $G$, and $R$ the
commutative ring $\Z/ n\Z$.  Informally, this construction is obtained
by (i) replacing an abelian group by the tensor product with $R$; (ii)
replacing the collection of abelian groups given by the filtration
quotients of $G$ obtained from the descending central series for $G$
by their tensor product with $R$; and (iii) reassembling this data
into a single group $G \otimes_{\mathbb Z} R$ via characteristic
classes of central extensions.  Details are given below.

Let $G$ denote a discrete group with $\Gamma^n$ the $n$-th stage of
the descending central series for $G$, defined inductively by
$\Gamma^1=G$, and $\Gamma^{n+1}=[\Gamma^n,G]$ for $n\ge 1$.  Since the
descending central series quotient $\Gamma^n/ \Gamma^{n+1}$ is an
abelian group, it is a module over the integers, and $\Gamma^n/
\Gamma^{n+1} \otimes_{\mathbb Z} R$ is defined for any commutative
ring $R$.  Observe that there is central extension
\[
1 \to \Gamma^n/ \Gamma^{n+1} \to G/\Gamma^{n+1} \to
G/\Gamma^{n} \to 1.
\]

Now define $ G/\Gamma^2 \otimes_{\mathbb Z} R = \Gamma^1/ \Gamma^{2} 
\otimes_{\mathbb Z} R$, and assume that 
$G/\Gamma^n \to\  G/\Gamma^n \otimes_{\mathbb Z} R$ has been
defined along with a commutative diagram
\[
\begin{CD}
K( G/\Gamma^n, 1) @>>>  K(G/\Gamma^n \otimes_{\mathbb Z} R,1)  \\
  @VVV                            @VVV \\
K(\Gamma^n/ \Gamma^{n+1},2)       @>>>    K(\Gamma^n/ \Gamma^{n+1} 
\otimes_{\mathbb Z} R,2)
\end{CD}
\] 

Notice that the homotopy theoretic fibre of
\[
K(G/\Gamma^n \otimes_{\mathbb Z} R,1) \to 
K(\Gamma^n/ \Gamma^{n+1} \otimes_{\mathbb Z}R,2)
\] 
is a $K(H,1)$ space, where $H= G/ \Gamma^{n+1} \otimes_{\mathbb Z}R$. 
Denote this fibre by $K( G/ \Gamma^{n+1} \otimes_{\mathbb Z}R, 1)$. 
With some mild hypotheses, there is a homotopy commutative diagram
\[
\begin{CD}
K( G/\Gamma^{n+1}, 1) @>{}>>  K(G/\Gamma^{n+1} \otimes_{\mathbb Z} R,1)  \\
  @VVV                            @VVV \\
K(\Gamma^{n+1}/ \Gamma^{n+2},2) @>>{}>    
K(\Gamma^{n+1}/ \Gamma^{n+2} \otimes_{\mathbb Z} R,2)
\end{CD}
\] 

Then define $K(G \otimes_{\mathbb Z} R,1)$ as the inverse limit of the
tower of the $K(G/\Gamma^n \otimes_{\mathbb Z} R,1)$.  The reader is
cautioned that this completion is different (on the level of
fundamental groups) from the completion usually used in the theory of
discrete groups.  This is the tensor product as defined in Bousfield
and Kan \cite{BK} for nilpotent groups.

Returning to the pure braid group, consider the abelianization map
$P_n \to \Z^{m}$, where $m=\binom{n}{2}$ as above.  There are maps
$P_n \to \Z^3$ constructed from the cohomology algebra of $ P_n $ as
follows.  Recall that the cohomology of $P_n$ is the quotient of the
exterior algebra generated by one-dimensional classes $A_{i,j}$ for $
1 \leq j < i \leq n$ by the ideal generated by $A_{i,j}A_{i,t} -
A_{t,j}\cdot A_{i,t} + A_{t,j}\cdot A_{i,j}$ for $1\le j<t<i \le n$.

Write $I(i,t,j)$ for any sequence of integers $1 \leq j<t<i \leq n$.
Define $p_{I(i,t,j)}: P_n \to \Z^{3}$ by the associated cohomology
class.  That is, the cohomology of $\mathbb Z^{3}$ has three
fundamental cycles $\iota_s $ for $s=1,2,3$ which are of degree $1$.
Define $p_{I(i,t,j)}$ by the formula
\[
p_{I(i,t,j)}^*(\iota_s)=
\begin{cases}
A_{i,j} & \text{if $s=1$,}\\
A_{i,t} & \text{if $s=2$,}\\
A_{t,j} & \text{if $s=3$.}
\end{cases}
\]

An associated homomorphism $\rho_{I(i,t,j)}: P_n \to \pi_3$ is
constructed as follows.  Define $\pi_{I(i,t,j)}: P_n \to \Z^{6}$ to be
the composite
$\pi_{I(i,t,j)} = \sigma \circ \Delta \circ p_{I(i,t,j)}$, where
$\sigma:\Z^6 \to \Z^6$ is the permutation given by the formula
$\sigma(n_1,n_2,n_3,n_4,n_5,n_6) = (n_1,n_3,n_2,n_5,n_4,n_6)$, and
$\Delta:\Z^3 \to \Z^6$ is defined by the formula
$\Delta(n_1,n_2,n_3) = (n_1,-n_1, n_2, n_2,n_3,n_3)$.

\begin{lem}
The characteristic class of the extension for the generalized
Heisenberg group $\pi_3$ pulls back to zero along $\pi_{I(i,t,j)}^*$,
that is $\pi_{I(i,t,j)}^*(\chi_3) = 0$.  Thus there is a lift to
$\rho_{I(i,t,j)}: P_n \to \pi_3$.  The product
\[
\rho =
\prod_{1 \le j < t < i \le n}\rho_{I(i,t,j)}: P_n
\to (\pi_3)^{\binom{n}{3}}
\]
of the $\binom{n}{3}$ maps $\rho_{I(i,t,j)}: P_n \to \pi_3$ induces a
surjection in cohomology.
\end{lem}
\begin{proof}
Notice that $\pi_{I(i,t,j)}^*(\chi_3) = A_{i,j}\cdot A_{i,t} - A_{t,j}
\cdot A_{i,t} + A_{t,j}\cdot A_{i,j}$ for $ j<t<i$.  This is precisely
the three-term relation for the cohomology algebra of $P_n$.  In
addition, the elements $A_{i,j}$ are in the image of $\rho^*$.  The
lemma follows.
\end{proof}

There are natural Spin representations obtained from the maps $\rho$
as follows.  Consider the mod $2$ reductions of the maps
$\pi_{I(i,t,j)}$ to obtain 
$\bar \pi_{I(i,t,j)}: P_n \to (\Z/ 2\Z)^{6}$.  The group $(\Z/2\Z)^6$
is a maximal elementary abelian $2$-group in $SO(7)$ while the
characteristic class defining $Spin(7)$ pulls back to that defining
$\pi_3$.  Thus, there are embeddings of 
$\pi_3 \otimes \Z/2\Z = \bar \pi_3 $ in the group $Spin(7)$, where
the tensor product is as defined above.

There are maps given by taking products of the above composites:
\[
\psi_{I(i,t,j)}: P_n \to Spin(7)
\qquad \text{and} \qquad
\Psi: P_n \to (Spin(7))^{\binom n3} \to Spin( 7\cdot{\binom n3}).
\]
Notice that the maps $\psi_{I(i,t,j)}$ are non-trivial.  Furthermore,
these maps with targets given by the Spinor groups do not factor
through the abelianization of $P_n$.

\begin{rem} Taking products, there are homomorphisms to products of
Heisenberg groups that factor through $\Gamma^1(P_n) / \Gamma^3(P_n)$
associated to the descending central series for $P_n$.  This
homomorphism with source $\Gamma^1(P_n) / \Gamma^3(P_n)\otimes \Z/2\Z$
is a monomorphism by inspection in case $n=3$.  For $n>3$, 
verification that this map is a monomorphism is left to the reader.
\end{rem}

The above discussion is summarized in the following.

\begin{thm} \

\begin{enumerate}
\item There are non-trivial Spin representations of $P_n$ given by
\[
\psi_{I(i,t,j)}: P_n \to Spin(7)
\qquad \text{and} \qquad
\Psi: P_n \to {Spin(7)}^{\binom n3} \to Spin( 7\cdot{\binom n3}).
\]
Furthermore, there is
a map, which is an embedding for $n= 3$,

\[
{\Gamma^1(P_n) / \Gamma^3(P_n) \otimes \Z/2\Z }
\to Spin(7)^{\binom n3}.
\]

\item Every Spin representation of $\Gamma$ where $\Gamma$ is
homologically toroidal induces the trivial map in real $K$-theory.
\end{enumerate}
\end{thm}

Notice that this last theorem gives many choices of non-trivial Spin
representations of the pure braid group for which the corresponding
fibre bundles are trivial.

\section{The Burau representation, $K$-theory, and the stable braid
group} \label{sec:FullBraidGroup}

Recall that Proposition \ref{prop:FullBraid} concerns the Burau
representation of the full braid group.

Let $\sigma_{k}$, $1\le k \le n-1$, denote the standard generators of
the full braid group $B_{n}$.  The Burau representation is the
homomorphism
\[
b: B_n \to\ GL(n, \mathbb Z[t,t^{-1}])
\]
given on generators by
\[
b(\sigma_{k}) = 
\begin{pmatrix}
\II_{k-1} & 0   & 0 & 0 \\
0         & 1-t & t & 0 \\
0         & 1   & 0 & 0 \\
0         & 0   & 0 & \II_{n-k-1}
\end{pmatrix},
\]
where $\II_{m}$ denotes the $m\times m$ identity matrix.  Among other
interpretations, this representation may be realized as the action of
$B_{n}$ on the one-chains of an infinite cyclic cover of a bouquet of
circles induced by the Artin representation.  See Birman \cite{Bi} for
a detailed account of the Burau representation.

To prove Proposition \ref{prop:FullBraid}, first consider the
specialization at $t=1$ of the Burau representation.  Notice that this
specialization is given by the natural map of the braid group to the
symmetric group followed by the natural inclusion of the symmetric
group in $GL(n,\Z)$.

Complexification of this representation gives the trivial bundle via
the Vandermonde matrix.  A complicated proof of this fact was
simplified in \cite{CMM}, and is included here for convenience of the
reader.

This bundle is given by the complex $n$-plane bundle
\[
F(\C,n)\times_{\Sigma_n}\C^n  \to  F(\C,n)/{\Sigma_n}.
\]
A trivialization is as follows.  Let $z = (z_1,\dots,z_n)$ be a point
in configuration space and $x = (x_1,x_2,\dots,x_n)$ in $\C^n$.
Define $\lambda: F(\C,n) \times \C^n \to F(\C,n) \times \C^n $ by the
formula $\lambda(z,x) = (z,y)$ where $y = (y_1,y_2,\dots,y_n)$ is
given by $y_i = \sum_{j=1}^n(z_i)^{j-1}x_i.$ Thus there is an
induced bundle isomorphism
\[
F(\C,n) \times_{\Sigma_n}\C^n \to F(\C,n)/{\Sigma_n} \times \C^n.
\]
The bundle on the right-hand side is trivial.  Thus the bundle on the
left-hand side is also trivial, as is the tensor product with any
other bundle.

Since $\C^*$ is path-connected, any specialization of the Burau
representation induces a map in $K$-theory analogous to that induced
by setting $t = 1$.  Thus parts (1) and (3) of Proposition
\ref{prop:FullBraid} follow.

To prove part (2), recall that Formanek \cite{F} shows that any
irreducible representation in $\Irr( B_n, GL(m,\C))$ for $m < n$, and
$ n > 6 $, is given by a specialization of a tensor product of a
one-dimensional representation with the reduced Burau representation.
The proposition follows.

Next consider a finite dimensional unitary representation of the
stable braid group
\[
\rho: B_{\infty} \to U(n).
\]
Notice that any invariant subspace of $\mathbb C^n$ admits an
orthogonal complement as the representation is unitary.  Thus a
unitary representation splits as a sum of irreducible representations. 
Assume that $m \gg n $.  In \cite[Lemma 9]{F}, Formanek shows that an
irreducible complex representation of $B_m$ of dimension $m-1$ for $m
> 2$ does not extend to $B_{m+2}$.  Thus any finite dimensional
irreducible unitary representation of $B_{\infty}$ of dimension at
least two is trivial.  Consequently, these representations give maps
which are trivial after passage to classifying spaces, and they induce
trivial maps on $K$-theory.

The remaining case is $ n = 1 $.  Since $U(1)$ is abelian, any
one-dimensional unitary representation $\rho: B_{\infty} \to U(1)$
factors through $\Z$, the abelianization of $B_{\infty}$, and thus
gives a trivial element in complex $K$-theory.  The next proposition
follows at once.

\begin{prop}
Let $\rho: B_{\infty} \to U(n)$ be a homomorphism.  The induced
element in complex $K$-theory $B\rho: B_{\infty} \to BU$ is trivial.
\end{prop}

The next lemma is included for future use as it provides a method for
considering general maps from $BB_{\infty}$ to $BU$.

\begin{lem}
The natural map from the direct limit $\lim BB_n \to BB_{\infty}$
induces a map
\[
[BB_{\infty},X] \to \varprojlim_n [BB_n,X]
\]
which is an isomorphism of groups whenever $X$ is an infinite loop
space.
\end{lem}
\begin{proof}
Since $X$ is an infinite loop space, there is a homotopy equivalence
$\Omega^q X_q \to X$ for some choice of space $X_q$, and any strictly
positive integer $q$.  Thus the induced map $[BB_{\infty},X] \to
[\Sigma^q BB_{\infty},\Omega X_{q+1}]$ is an isomorphism of groups.
Recall \cite{C2} that $\Sigma^{2L} BB_{\infty}$ is homotopy equivalent
to $\Sigma^{2L} BB_{L} \vee \Sigma^{2L}Y_L$ where $Y_L$ is the cofibre
of the natural map $BB_{L} \to BB_{\infty}$, and that $Y_L$ is
$\lfloor L/2 \rfloor$-connected.

Recall Milnor's $lim^1$ exact sequence from \cite{Mil}, as described
in Bousfield-Kan \cite{BK} on the level of function spaces: let $E =
colim_n E_n$ be a filtered space, where $(E_n, E_{n-1})$ is an NDR
pair in the terminology of Steenrod; then there is an exact sequence
of groups
\[
1 \to lim^1[E_n,X] \to [E, X] \to \varprojlim_n [E_n,X] \to 1.
\]
Next let $E_n$ denote $BB_n$.  By the previous remarks concerning the
braid groups, the restriction maps $[E_n,X] \to [E_{n-1},X]$ are split
epimorphisms of groups (as $X$ is assumed to be an infinite loop
space).  The Mittag-Leffler condition is satisfied, and thus Milnor's
$lim^1$ term is zero (see \cite[Prop. 2.4, Cor. 3.3]{BK}).  The
lemma follows.
\end{proof}

\section{An analogue of a classifying space for $KO^{0}_{\rep}$}
\label{sec:ClassifyingSpace}
Let $\Gamma$ be any group, with $ g: \Gamma \to O(n)$ a group
homomorphism.  Composing $g$ with the natural map $O(n) \to O$ yields
an induced map $G:\Gamma \to O$.  There is a function
\[
\Phi:[B\Gamma,BO] \to H^1(B\Gamma;\Z/2\Z) \oplus H^2(B\Gamma;\Z/2\Z)
\]
obtained by evaluating a map on $w_1$ and $w_2$, the first two
Steifel-Whitney classes.

The purpose of this section is to give a natural description of the
above function in terms of spaces.  Recall from real Bott periodicity
that there is a homotopy equivalence $BO \to \Omega(SU/SO)$, and that
$BSpin$ is the $2$-connected cover of $BO$.

Consider the map $Sq^2:K(\Z/2\Z,2) \to K(\Z/2\Z,4)$ given by sending
the fundamental cycle in the mod-$2$ cohomology of $K(\Z/2\Z,4)$ to
the cup square of the fundamental cycle for the mod-$2$ cohomology of
$K(\Z/2\Z,2)$.  Let $E$ denote the homotopy theoretic fibre of this
map.  The space $E$ is the first two stages of the Postnikov tower for
the single delooping of $BO$.

Notice that the loop space of $E$ is a product given by $K(\Z/2\Z,1)
\times K(\Z/2\Z,2)$, but this decomposition does not preserve the loop
structure.  Thus the sets $[X,\Omega(E)]$, and $[X,K(\Z/2\Z,1) \times
K(\Z/2\Z,2)]$ are isomorphic, but may have different group structures.

A standard and classical exercise carried out below gives a loop map
$\text{ev}: SU/SO \to E$ which induces an isomorphism on the first
five homotopy groups.  Identifying $BO$ with $\Omega (SU/SO)$, there
is a double loop map $\Omega(\text{ev}):BO \to \Omega(E)$.  The reader
is referred to Cartan's computations in \cite{Car} for the details of
the proof.

There are isomorphisms
\[
H^*(U/O; \mathbb F_2) \cong E[x_1,x_2, \dots,x_n,\dots]
\quad \text{and} \quad
H_*(U/O; \mathbb F_2) \cong \mathbb F_2[y_1,y_3,\dots,y_{2n+1},\dots].
\]
Thus the unique non-trivial map $SU/O \to K(\Z/2\Z,2)$ lifts to $E$,
and induces an isomorphism on the first non-vanishing homotopy group.
Furthermore, all maps are loop maps as $Sq^2$ is stable.
Consequently, there is a lift which is a loop-map, and induces an
isomorphism up to $H_3$.  This map induces an isomorphism on the first
three homotopy groups.  This suffices by Bott periodicity and the
definition, as the next possibly non-vanishing group is in dimension
five.

\begin{prop} \

\begin{enumerate}
\item The homotopy theoretic fibre of $BO \to \Omega(E)$ is $BSpin$.
Thus for any pointed space $X$, there is a long exact sequence of
abelian groups
\[
\cdots \to [X,\Omega^2(E)] \to [X,BSpin] \to [X,BO] \to [X,\Omega(E)]
\to [X,B^2Spin] \to \cdots.
\]

\item For any path-connected $CW$-complex $X$ there is a short exact
sequence of abelian groups
\[
0 \to H^2(X;\Z/2\Z) \to [X,\Omega(E)] \to H^1(X;\Z/2\Z) \to 0
\]
which is split as sets.  Thus $H^2(X;\Z/2\Z)$ acts on $[X,\Omega(E)]$
which is isomorphic to $H^1(X;\Z/2\Z) \times H^2(X;\Z/2\Z)$ as an
$H^2(X;\Z/2\Z)$-set.

\item In case $X=B\Gamma$ for a homologically toroidal group $\Gamma$,
the composite
\[
KO^0_{\rep}(B\Gamma)\to KO^0(B\Gamma) \to [B\Gamma,\Omega(E)]
\]
is an isomorphism of groups, and $[B\Gamma,\Omega(E)]$ is isomorphic
to $H^1(B\Gamma;\Z/2\Z) \oplus H^2(B\Gamma;\Z/2\Z)$ with this choice
of isomorphism given by evaluating a map on $w_1$ and $w_2$, the first
two Stiefel-Whitney classes.
\end{enumerate}
\end{prop}
\begin{proof}
Most of the proof has been outlined above.  Notice that the natural
map $\Omega(\text{ev}):BO \to \Omega(E)$ induces an isomorphism on the
first three homotopy groups.  Since $E$ has two non-trivial homotopy
groups, and the map $\Omega(\text{ev})$ is an isomorphism of these two
homotopy groups, the map induces an isomorphism on the first $3$
homotopy groups by Bott periodicty.  Thus the fibre of
$\Omega(\text{ev}):BO \to \Omega(E)$ is $BSpin$, and the long exact
sequence of part (1) follows at once.

To prove part (2), notice that there is a fibration $K(\Z/2\Z,3) \to E
\to K(\Z/2\Z,2)$ which admits a section after looping.  Thus there is
a short exact sequence
\[
0 \to H^2(X;\Z/2\Z) \to [X,\Omega(E)] \to H^1(X;\Z/2\Z) \to 0.
\]

Notice that $ KO^0_{\rep}(B\Gamma)$ is the subgroup of $KO^0(B\Gamma)$
generated by elements for which at least one of the first two
Stiefel-Whitney classes is non-zero.  This is precisely the group
$[X,\Omega(E)]$.  Since every element in $ KO^0_{\rep}(B\Gamma)$ has
order $2$ as shown in Proposition \ref{prop:CtrivialPROP}, this
sequence is split as groups, part (3) follows.
\end{proof}

\section{Final remarks: representations into $GL(n,R)$}
Rather than restricting to either orthogonal or unitary
representations of a group, it is natural to also consider more
general representations with values in $GL(n,R)$ for $R$ either the
real, or complex numbers, and to pose the same questions about
$K$-theory.  Some remarks concerning more general represntations are
included in this final section.  With this change, the situation
regarding $K$-theory, and representations might be more complicated
than that given by orthogonal or unitary representations depending on
the choice of group.  This will be illustrated in three ways below.

As a first example, consider representations of a surface group $\pi$
in $SL(2,\mathbb R)$ so that the action of $\pi$ on the homogeneous
space $SL(2,\mathbb R)/SO(2)$, regarded as the upper half-plane, is
properly discontinuous.  The existence of such representations dates
back to the nineteenth century.  Notice that such a represention
cannot compress to $SO(2)$ as the associated action of a purported
compression is not properly discontinuous.

A second type of situation occurs from the fact that an arbitrary
representation in $GL(n, \mathbb C )$ may not be a sum of irreducible
representations if the group $\pi$ is not finite.  In the arguments in
Section \ref{sec:circles} above, an elementary but essential use is
made of the fact that unitary representations split as sums of
irreducible representations.

Next, consider the $6$-stranded braid group for the $2$-sphere denoted
by $\Gamma^6$.  This group is a quotient of the mapping class group
$\Gamma_{2}$ for a genus $2$ surface.  The natural homomorphism
$\Gamma_2 \to\ Sp(4,\mathbb Z)$ descends to a homomorphism
\[
\phi: \Gamma^6 \to\ PGL(2, \mathbb C ),
\]
obtained by regarding $Sp(4,\mathbb Z)$ as a subgroup of 
$Sp(4,\mathbb R)$ which, in turn, is a subgroup of $GL(2, \mathbb C)$. 
A number of properties of this representation are listed below.  Of 
these, (1) is easily verified, while (2) may be established using
(1) and the results of \cite{BC}.

\begin{enumerate}
\item The map $\phi$ induces a monomorphism
in mod-$2$ cohomology on the level of classifying
spaces $B\phi: B\Gamma^6 \to\ BPGL(2, \mathbb C)$
by restricting to an elementary abelian $2$-group
of rank $2$ in $\Gamma^6$.
\item The map $B\phi$ is not homotopic to a map
$B\rho$ arising from a representation
$$\rho:\Gamma^6 \to\ SO(3).$$
\item The group $SO(3)$ is a maximal compact
subgroup of $PGL(2, \mathbb C)$.
\item The representation $\phi$ is thus not orthogonal,
but is both natural and interesting.
\item The pure braid group associated to $\Gamma^6$
is homologically toroidal, and the results proven here
for orthogonal, and unitary representations applies
to this group.
\end{enumerate}

The point of these remarks is that there are representations of the
braid groups which are not orthogonal and yet contain relevant
topological information.  While some of methods of this paper apply to
these more general representations in the case of pure braid groups, a
comprehensive analysis of the $K$-theoretic contributions of these
representations remains to be carried out.

\begin{ack}
Some of this work was carried out while the second two authors visited
the Department of Mathematics at the University of Wisconsin-Madison
in the Spring of 1999.  These authors thank the Department and the
University for their hospitality, and for providing an exciting and
productive mathematical environment.
\end{ack}

\bibliographystyle{amsalpha}

\end{document}